\documentclass[10pt, oneside]{article}   	% use "amsart" instead of "article" for AMSLaTeX format
\usepackage{geometry}                		% See geometry.pdf to learn the layout options. There are lots.
\usepackage{graphicx}				% Use pdf, png, jpg, or eps§ with pdflatex; use eps in DVI mode
\usepackage{amssymb}

\title{Remarks on the nonvanishing of cohomology groups for perverse sheaves on abelian varieties}
\author{Rainer Weissauer}
\date{}							% Activate to display a given date or no date

\begin{document}
\maketitle
%\section{}
%\subsection{}

\bigskip\noindent

Let $X$ be an abelian variety over an algebraically closed field $k$
of dimension $g$ and let $K$ be an irreducible perverse sheaf 
in $D_c^b(X,\Lambda)$ for $\Lambda={\overline \mathbb Q}_\ell$.
If the base field $k$ has positive characteristic, we assume that $K$ is defined over
a field that is finitely generated over its prime field with $\ell$ different from the characteristic.
Suppose that not all cohomology groups 
$H^\nu(X,K)$ are zero and let denote $d(K)= \max\{ \nu\ \vert \ H^\nu(X,K)\neq 0\}$. Notice $d(K)\geq 0$, by the Hard Lefschetz Theorem. 

\medskip
{\bf Theorem}. {\it  For $d=d(K) > 0$  we have 
$$ \dim_\Lambda(H^{d-1}(X,K)) \ > \ 2d/(d+g) \cdot \dim_\Lambda(H^d(X,K)) \ .$$ 
If furthermore $X$ is a simple abelian variety, then $\dim_\Lambda(H^{d-1}(X,K)) > d \cdot \dim_\Lambda(H^d(X,K))$.}

\bigskip
{\bf Remark}. By the Hard Lefschetz Theorem an immediate consequence of this theorem is the assertion: $H^\nu(X,K)\!\neq\! 0$ if and only if $\nu\! \in\! [-d(K),d(K)]$. So for character twists $K_\chi$ [KrW] the sets $V_i(K)\! =\! \{\chi \vert H^i(X,K_\chi)\neq 0\}$  satisfy $V_{i+1}(K) \subseteq V_{i}(K)$ for all $i\geq 0$. For an arbitrary
projective smooth variety $Y$ over $k$ with Albanese morphism $f:Y\to X$
and a perverse sheaf $L$ on $Y$ the decomposition theorem gives $Rf_*(L)\cong
\bigoplus {}^p\! H^i(Rf_*(L))[-i]$ and $H^\nu(Y,L)\cong \bigoplus_{j+i=\nu} H^j(X, {}^p\! H^i(Rf_*(L))$.
From the relative Hard Lefschetz Theorem and the theorem above applied to the irreducible constituents $K$ of the semisimple perverse 
cohomology sheaves ${}^p\! H^i(Rf_*(L))$ we therefore obtain 

\medskip
{\bf Corollary 1}. {\it Let $L$ be an irreducible perverse sheaf $L$ on a smooth projective variety $Y$ with $d= \max\{ \nu\ \vert \ H^\nu(Y,L)\neq 0\}$. 
Suppose the Albanese morphism $f\!:\! Y\!\to\! X$ is not trivial and suppose   
$H^d(Y,L) \neq H^0(X,{}^p\! H^d(Rf_*(L))) $ (e.g. this is the case if the fibers of $f$ have dimension $<d$).
Then $H^\nu(Y,L)\neq 0$ holds if and only if $\nu \in [-d,d]$.}  

\medskip
{\it Proof of the theorem}. First suppose that $K$ is negligible, i.e. of the form $K \cong \pi^*(Q)[q]$
for a perverse sheaf $Q$ on a quotient abelian variety $\pi: X \to X/A$ 
defined by an abelian subvariety $A \subseteq X$ of dimension $q>0$.
Then $d=d(K) =d(Q) + q$ since $H^\bullet(X,K)\cong \bigoplus_{i=0}^{2q} {2q \choose i} \cdot
H^\bullet(X/A,Q[i+q])$. Hence $H^{d}(X,K) \cong H^{d(Q)}(X/A,Q)$ and
$H^{d-1}(X,K) \cong 2q \cdot  H^{d(Q)}(X/A,Q) \oplus H^{d(Q)-1}(X/A,Q)$.
Since $2q > 2d/(d+g)$, our claim follows in this case; similarly $2q=2g >d$
in the case where $X=A$ is simple. Therefore we now make the

\medskip
{\it Assumption}. {\it Suppose $K$ is irreducible, but not negligible.
Furthermore suppose $d>0$}.

\medskip
For the perverse sheaf $K$ on $X$ consider  the Laurent polynomial
$h_t(X,K) = \sum_\nu a_\nu t^\nu$ defined by 
$a_\nu = \dim_\Lambda(H^\nu(X,K))$. Then $d=d(K)$ is the largest integer 
$\nu$ such that $a_\nu\neq 0$.

\medskip
Choose an integer $r$ minimal such that $r \cdot d > g$. Hence $r >1$
and $r\cdot d < g + d$. The $r$-th convolution power of $K$
is a direct sum of a perverse sheaf $K_r$ on $X$ and a finite direct sum
of complexes $L_\mu[n_\mu]$ with negligible perverse sheaves $L_\mu$
on $X$ of the form:
\begin{itemize}
\item $L_\mu = \pi_\mu^*(Q_\mu)[g_\mu]$ for irreducible not negligible perverse sheaves $Q_\mu$ on
$X/A_\mu$ 
\item $\pi_\mu\!: \! X \to X/A_\mu$ is the quotient by an abelian subvariety 
$A_\mu$ of $X$ of dimension $g_\mu >0$.
\end{itemize} 
This follows from [KrW], [W], and for this assertion we have
to assume that the perverse sheaf $K$ is defined over a finitely generated field over the prime field in the case  
of positive characteristic [W]. 

\medskip
Then $h_t(X,L_\mu[n_\mu])= \sum_\nu \dim(H^\nu(X,L_\mu[n_\mu])\cdot t^\nu 
= \sum_{\nu\leq d_\mu} b_{\mu \nu} t^\nu$ for integers $b_{\mu\nu}\geq 0$, and
we may assume $b_\mu = b_{\mu d_\mu} \geq 1$ since we can ignore cohomologically trivial summands in the following. Let $T$ denote 
the  set of all indices $\mu$ such that $d_\mu + g_\mu = r \cdot d$ holds.
By well known cohomological bounds [BBD], the cohomology of an irreducible perverse sheaf on $X$  vanishes in degrees $\geq g$
unless it is negligible. Since $r \cdot d \geq g$, the K\"unneth formula in the form $H^\bullet(X,K^{*r}) 
\cong H^\bullet(X,K)^{\otimes r}$ and a 
comparison
of coefficients at $t^{rd}$   implies 
$$   (a_d)^r = \sum_{\mu \in T} \  b_\mu  \ .$$ 
Similarly, now using $r\cdot d  \geq g +1$, by comparing coefficients
at $t^{rd-1}$ we obtain 
$$ r \cdot a_{d-1} (a_d)^{r-1} \ \geq \ \sum_{\mu \in T} 2 g_\mu b_\mu \ \geq\  2 \cdot \min_\mu(g_\mu) \cdot (a_d)^r   \  .   $$
Indeed, the second equality follows from the formula $\sum_{\mu\in T} b_\mu = (a_d)^r$ above. 
For the first inequality we exploited the fact that all coefficients $b_{\mu\nu}$ in
$h_t(X,L_\mu[n_\mu]) = (t+2+t^{-1})^{g_\mu} \cdot h_t(X/A_\mu, Q_\mu[n_\mu])
= (t^{g_\mu} + 2g_\mu t^{g_\mu - 1} + \cdots )(b_\mu t^{d_\mu} + \cdots )$  
are nonnegative.
We conclude
$$       a_{d-1} \ \geq \ \frac{ 2 \min_\mu(g_\mu) }{r}  \cdot a_d \ \geq \ \frac{2}{r}\cdot  a_d \ > \ \frac{2d}{g+d}\cdot a_d $$
where the last inequality follows from $r\cdot d < g + d$. If $X$ is simple, then
$\min_\mu(g_\mu)  = g$ and hence $a_{d-1} \geq \frac{2g}{r} a_d$. Now $r\cdot d < g + d < 2g$
implies
$    a_{d-1} >  d \cdot a_d  $. QED

\medskip
{\bf Remark}. $d(K)$ for the intersection cohomology sheaf $K$
of an irreducible subvariety $Y$ of $X$ is the dimension of $Y$.
In this case there exist stronger geometric estimates 
than those from the theorem above. 
However, already when $Y$ is a variety of maximal Albanese
dimension and $K$ is an arbitrary irreducible constituent of the direct image 
of the intersection cohomology perverse sheaf on $Y$ under the Albanese morphism $f\!:  Y \to X\!=\! Alb(Y)$ I am not aware of estimates of the above form in the literature.

\medskip
Next, consider a finite Galois morphism $$\pi: \tilde Y \to Y$$ between smooth complex varieties of dimension $n$ with Galois group $\Gamma$, where we view $\Gamma$ to act on $\tilde Y$  from the right. For every isomorphism class $\phi$ of irreducible representations $V_\phi$ of $\Gamma$ 
let $m_{\nu}(\phi)$ denote the multiplicity of 
the irreducible representation $\phi$ of $\Gamma$ on $H^{\nu+n}({\tilde Y},\mathbb C)$. 

\medskip
For simplicity, from now on suppose that $Y$ is projective and 
$f: Y \to Alb(Y)=X$ is a closed embedding.  Then 
the theorem above implies 

\medskip
{\bf Corollary 2}. {\it  If $d = d(K_\phi) >0$, 
then $m_{d-1}(\phi) > 2d m_d(\phi)/(d+g)  > 0$.} 

\medskip
{\it Proof}. For every class $\phi$  there exists an irreducible perverse sheaf $K_\phi$ on $Y$ and a $\Gamma$-equivariant isomorphism $H^{\bullet+n}(\tilde Y,\mathbb C) \cong \bigoplus_\phi V_\phi \otimes_{\mathbb C} H^\bullet(Y, K_\phi)$, where $\Gamma$ acts on $V_\phi$ by $\phi$ and trivially on 
$H^\bullet(Y, K_\phi)$. For unramified $\pi$, this immediately follows from [KiW], remark 15.3 (d).  
Applying this remark for the restriction of $\pi$ to $\pi^{-1}(U)$, for the open dense subset $U\subseteq Y$ obtained by removing the ramification divisor of $\pi$, by perverse analytic continuation in general it suffices to observe that for $\delta_{\tilde Y} = \mathbb C_{\tilde Y}[n]$
the semisimple perverse sheaf $\pi_*(\delta_{\tilde Y})$ on $Y$ has 
irreducible perverse constituents $K$ whose restriction to $U$ are nontrivial. To show this notice that $Hom(\pi_*(\delta_{\tilde Y}),K) = Hom(\delta_{\tilde Y}, \pi^!(K))$ vanishes if $K$ (and hence $\pi^!(K)$) is a perverse sheaf with support of dimension $< \dim(Y)$. Indeed, since $\delta_{\tilde Y}$ is an irreducible perverse sheaf with support of dimension $\dim(Y)$, $Hom(\delta_{\tilde Y}, \pi^!(K))$ is zero.
This being said, we obtain $m_\nu(\phi)=\dim(H^{\nu}(Y,K_\phi))$. Since $f$ is a closed immersion, the direct images of $K_\phi$ under the Albanese morphism again are irreducible perverse sheaves. So we can apply the theorem. QED

\medskip
Still suppose $\pi: \tilde Y \to Y$ is a Galois covering and $f: Y \to Alb(Y)$ is a closed embedding.
Since $\chi(Y,K_\phi) = \sum_\nu (-1)^\nu \dim H^\nu(Y,K_\phi)$,
for $\gamma\in \Gamma$ the trace
$ tr(\gamma)\ =\ \sum_\nu (-1)^\nu tr(\gamma; H^\nu(\tilde Y, \delta_{\tilde Y}))$
can be written 
$$ tr(\gamma) = \sum_\phi \ \chi(Y,K_\phi) \cdot tr(\gamma;V_\phi) \ .$$
$K_\phi$ has rank $dim(V_\phi)$ on $U$, thus generic rank $dim(V_\phi)$ on $Y$. 
Hence the characteristic cycle of the 
D-module on $Alb(Y)$ attached to $f_*(K_\phi)$ is a sum of irreducible Lagrangians cycles
containing the conormal Lagrangian cycle $\Lambda_{f(Y)} \subset T^*(X)$
with multiplicity $\dim(V_\phi)$.  
As shown in [FK], by the theorem of Dubson-Riemann-Roch this implies 
$\chi(K_\phi)= \chi(f_*(K_\phi)) \geq \dim(V_\phi) \cdot deg(\Lambda_{f(Y)})$.
Furthermore since $f(Y) \cong Y$ is smooth, the characteristic variety of $f(Y)$
is $\Lambda_{f(Y)}$ and hence 
$deg(\Lambda_{f(Y)})$  is the Euler-Poincare
characteristic $\chi_Y$ of the variety $Y$, again by [FK]. 
%This implies $\chi(K_\phi) \geq dim(V_\phi)\cdot \chi_Y$ (where the Euler characteristic $\chi_Y$ enters as the degree of the Gauss map and $dim(V_\phi)$
%as the generic rank that gives the multiplicity of the generic component of the characteristic cycle; notice that since $Y$ is smooth $deg(\Lambda_Y)=\chi_Y$), 
%since the other characteristic cycles only enlarge the value by the theorem of Dubson-Riemann-Roch. 
This implies
$$ tr(\gamma)  = 
\sum_\nu \ (\dim(V_\phi)\chi_Y + a_\phi) \cdot \phi(\gamma) \   $$ for certain integers $a_\phi\geq 0$.
Hence the virtual representation defined by $tr(\gamma)$ is $\chi_Y$ times the regular representation
of $\Gamma$ plus a true representation of $\Gamma$. Notice that $\chi_Y \geq 0$ holds by [FK]
and our assumptions on $f$. 

\medskip
{\bf Remark}. In the case of surfaces $Y$, for a nontrivial irreducible representation $\phi$ of $\Gamma$ from this Chevalley-Weil type trace formula  we obtain the estimate $m_0(\phi) - 2m_1(\phi)= \dim(V_\phi)\chi_Y + a_\phi \geq 0$. 
So, for surfaces and nontrivial $\phi$ under the assumptions before corollary 2, this improves the previous estimate $m_0(\phi) \geq 2m_1(\phi)/(g+1)$ of corollary 2.

\bigskip\noindent
{\bf References}:
 
\medskip
[BBD] Beilinson A., Bernstein J.,  Deligne P., {\it Faisceaux pervers}, Asterisque 100 (1982).

\medskip
[FK] Franecki J., Kapranov M., {\it The Gauss map and a noncompact Riemann-Roch formula for constructible sheaves on semiabelian varieties}, Duke Math. J. 104 no. 1 (2000) 171-180.

\smallskip
[KiW] Kiehl R., Weissauer R., {\it Weil conjectures, Perverse Sheaves and l-adic Fourier Transform},
Springer Verlag, Ergebnisse der Mathematik 42 (2001).

\smallskip
[KrW] Kr\"amer Th., Weissauer R., {\it Vanishing theorems for constructible sheaves on abelian varieties},
J. Alg. Geom. 24 (2015), 531 - 568.

\smallskip
[W] Weissauer R., {\it Vanishing theorems for constructible sheaves on abelian varieties over finite fields}, To appear in Math. Annalen.

\goodbreak

\end{document}